\documentclass[11pt,a4paper, reqno]{amsart}
\usepackage{amssymb, paralist}
\usepackage{amsfonts, amsmath, wasysym}
\usepackage{latexsym}
\usepackage{graphicx, color}
\usepackage{indentfirst}
\usepackage{stmaryrd}

\usepackage{amssymb}
\usepackage{amsmath}
\usepackage{graphicx, color}
\usepackage{indentfirst}

\newcommand{\N}{\mathbb{N}}

\newcommand{\re}{\mathbb R}

\newcommand{\<}{\left<}
\newcommand{\fracc}{\displaystyle \frac}
\newcommand{\ds}{\displaystyle}
\renewcommand{\>}{\right>}
\newcommand{\flecha}{\longrightarrow}

\newcommand{\mc}[1]{\mathcal{#1}}

\def\tr{{\rm tr}}

\def\fle{\rightarrow}
\def\parcial#1#2{\fracc{\partial #1}{\partial#2}}
\def\deri#1#2{\fracc{d #1}{d#2}}

\def\({\left (}
\def \){\right)}

\newtheorem{thm}{Theorem}[section]

\newtheorem*{thm*}{Theorem} 
\newtheorem{cor}[thm]{Corollary}\newtheorem{prop}[thm]{ Proposition}

\theoremstyle{definition}

\newtheorem{defin}[thm]{Definition}

\newtheorem*{rems*}{Remarks}
\theoremstyle{remark}
\newtheorem{rem}[thm]{Remark}

\newcommand{\eps}{\ensuremath{\varepsilon}}

\def\p{\varphi}

\newtheorem{teor}{\hspace{12pt} Theorem}

\newtheorem{lema}[teor]{\hspace{12pt} Lemma}

\numberwithin{teor}{section}

\newcommand{\be}{\begin{enumerate}}
\newcommand{\ee}{\end{enumerate}}
\newcommand{\bi}{\begin{itemize}}
\newcommand{\ei}{\end{itemize}}
\newcommand{\bd}{\begin{description}}
\newcommand{\ed}{\end{description}}
\newcommand{\bec}{\begin{equation}}
\newcommand{\eec}{\end{equation}}
\newcommand{\ba}{\begin{array}}
\newcommand{\ea}{\end{array}}
\newcommand{\bt}{\begin{thm}}
\newcommand{\et}{\end{thm}}
\newcommand{\bdem}{\begin{proof}}
\newcommand{\edem}{\end{proof}}
\newcommand{\bl}{\begin{lema}}
\newcommand{\el}{\end{lema}}
\newcommand{\bnp}{\begin{rem}}
\newcommand{\enp}{\end{rem}}
\newcommand{\bde}{\begin{defin}}
\newcommand{\ede}{\end{defin}}
\newcommand{\bnod}{\begin{rem}}
\newcommand{\enod}{\end{rem}}
\newcommand{\bp}{\begin{prop}}
\newcommand{\ep}{\end{prop}}
\newcommand{\bco}{\begin{cor}}
\newcommand{\eco}{\end{cor}}

\newcommand{\nn}{\nonumber}
\newcommand{\lb}{\label}

\newcommand{\ene}{\end{equation} }

\newcommand{\fint}{{-}\hspace*{-1.05em}\int}
\newcommand{\sfint}{{-}\hspace*{-0.90em}\int}

\def\l{\lambda}

\usepackage{bm}

\usepackage[draft]{optional}

\renewcommand{\(}{\left(}

\renewcommand{\>}{\right>}
\renewcommand{\)}{\right)}
\def\bal{\begin{align}}
\def\eal{\end{align}}

\numberwithin{equation}{section}
\def\be{\begin{equation}}
\def\ee{\end{equation}}

\def\oH{{\overline H}}

\def\tr{{\rm tr}}

\def\parcial#1#2{\frac{\partial #1}{\partial#2}}

\def\deri#1#2{\frac{d #1}{d#2}}

\def\flecha{\longrightarrow}
\def\fle{\rightarrow}

\def\ds{\displaystyle}

\def\tr{{\rm tr}}

\def\oH{{\overline H}}

\begin{document}




\title[Non-preserved curvature conditions under constrained mcf]{Non-preserved curvature conditions under constrained mean curvature flows}

\author{ Esther Cabezas-Rivas} 
\address{Institut f\"ur Mathematik, Goethe-Universit\"at Frankfurt 
\\ Robert-Mayer-Str.~10 \\ 60054 Frankfurt, Germany}
  \email{cabezas-rivas@math.uni-frankfurt.de}
  \author{ Vicente Miquel }
\address{Department of Geometry and Topology\\
  Fac. Mathematics. University of Valencia\\
  calle Dr. Moliner 50, 46100-Burjassot (Valencia) Spain}
  \email{vicente.f.miquel@uv.es}

\maketitle

\vspace{-1cm}
\begin{abstract}
We provide explicit examples which show that mean convexity (i.e.~positivity of the mean curvature) and positivity of the scalar curvature are non-preserved curvature conditions for hypersurfaces of the Euclidean space evolving under either the volume- or the  area preserving mean curvature flow. The relevance of our examples is that they disprove some statements of the previous literature, overshadow a widespread {\it folklore conjecture} about the behaviour of these flows and bring out the discouraging news that a traditional singularity analysis  is not possible for constrained versions of the mean curvature flow.
\end{abstract}

\bigskip



\section{Introduction and Main Results }\lb{In}


We consider a closed $n$-dimensional manifold $M$ smoothly embedded as a hypersurface of the Euclidean space $\re^{n + 1}$ by means of the parametrization $F_0: M \fle \re^{n + 1}$. We deform $M_0:= F_0(M)$ according to a {\it constrained Mean Curvature Flow} ({\sc mcf}), that is, we consider a family of isometric immersions $F: M \times [0, T) \fle \re^{n + 1}$ solving the initial value problem:
\begin{equation}\label{vpmcf}
\left\{\ba{ll}
\ds \parcial{F}{t} (p, t) = \big(h(t) - H(p, t)\big)\, N(p, t), & p\in M, \ t\in (0, T) \medskip \\ \quad F(\cdot, 0) =F_0 \ea \right.,
\end{equation}
where $N$ is the unit normal vector field pointing outward (if each evolving hypersurface $M_t= F(M,t)$ encloses a domain $\Omega_t$) and $H$ the mean curvature of $M_t$. 
More precisely, hereafter we study two types of constrained flows depending on the choice of the global term $h(t)$:

$\bullet$ The {\it Volume Preserving Mean Curvature Flow} ({\sc vp-mcf}) is obtained by defining $h(t)$ as the average mean curvature:
\bec \lb{oH}
h(t) = \oH=\ds\frac{ \int_{M} H d \mu_t }{|M_t| },
\eec
where $ d \mu_t$ is the induced Riemannian volume element of $M_t$,  and $|M_t| = \int_M d\mu_t$ denotes the corresponding $n$-volume (which we shall call {\it area}).    

$\bullet$ The {\it Area Preserving Mean Curvature Flow }({\sc ap-mcf}) is defined by choosing
\bec \lb{ap-mcf}
h(t) = \frac{\int_M H^2 \, d\mu_t}{\int_M H \, d\mu_t}.
\eec


The presence of the global term $h(t)$ in the flow equation \eqref{vpmcf} has two major consequences: 
\begin{enumerate}
\item for $h(t)$ defined as in \eqref{oH} the flow keeps
the enclosed volume constant while the area decreases as time evolves. On the other hand, if we choose $h(t)$ as in \eqref{ap-mcf}, any closed hypersurface moves preserving its area and increasing the volume enclosed by the hypersurface. 

\item  makes the usual techniques
in geometric flows (e.g.~the application of maximum principles) either fail or become more
subtle. 
\end{enumerate}

The resultant evolution problem (see Section \ref{pre} for a brief account of the previous literature) is particularly appealing -since from (1) it is specially
well suited for applications to the isoperimetric problem- and challenging because (2) causes
numerous extra complications; for instance, a basic property for the ordinary {\sc mcf} (the comparison principle) fails in general for \eqref{vpmcf}, e.g., an initially embedded
curve may develop self-intersections (cf.~\cite{MaSi}) under {\sc vp-mcf}. Hence the present knowledge of \eqref{vpmcf} is considerably poorer than that of the unconstrained evolution. This is even more evident when we talk about the analysis of singularities.

For the unconstrained {\sc mcf} a systematic study was started by G.~Huisken in \cite{Hu90}; in the latter and subsequent papers (e.g.~\cite{HuSi}) the preservation of the mean convexity ($H > 0$) played a key role. So a natural first step towards a similar analysis of singularities for the non-local flows would be to wonder if mean convexity is preserved or not under \eqref{vpmcf}. There is (even written) evidence that the community of experts was quite inclined towards believing in a positive answer:
\begin{enumerate}
\item[(a)] The preservation of $H > 0$  under {\sc vp-mcf} was claimed and {\it proved} in \cite[Lemma 3.11]{Wa}. 

\item[(b)] The corresponding statement and {\it proof} for the {\sc ap-mcf} is indeed published in \cite[Lemma 10]{DGG}.

\item[(c)] It is a folklore conjecture that {\it if the average mean curvature is bounded, we expect that the behaviour of singularities for the {\sc vp-mcf} after parabolic rescaling is the same as by {\sc mcf}} (quote adapted from \cite{Ath2}).
\end{enumerate}

Roughly speaking in this note we give counterexamples to  the {\it statements} in (a) and (b), that is, we find mean convex surfaces  in $\re^3$  which evolve under {\sc ap-} or {\sc vp-mcf} to surfaces with negative mean curvature at some points. Let us remark that, despite (a) and (b) study \eqref{vpmcf} in a certain non-euclidean ambient space, the ambient curvature plays no role in their {\it proof} of preservation of $H > 0$.

Furthermore, our examples of loss of mean convexity are all rotationally symmetric closed surfaces for which the average mean curvature is known to be bounded (cf. \cite[Lemma 8.4]{Atha-Kan}). Accordingly our examples also {\it throw cold water on} (c): it could still be true that the singularity behaviour is similar after a long while but for short times we prove that the evolution of mean convex hypersurfaces under \eqref{vpmcf} is quite different to that under ordinary {\sc mcf} already at the very beginning (that is, for short times).

In any case, our examples evidence that a singularity analysis following the traditional approaches (\cite{HuSi2}, \cite{Wh} and references therein) is not possible for the constrained versions of {\sc mcf}. The only remaining hope is to produce completely different techniques to face the study of singularities or to find another curvature condition (stronger than mean convexity but weaker than convexity) which is preserved under \eqref{vpmcf}. In the second direction, the next natural attempt is to wonder about preservation of the positivity of the scalar curvature ($R > 0$). We also address such a question in this paper (see the beginning of Section \ref{Smo} for reasons to justify why this is the next natural step) and show that this is again a condition that is not preserved under \eqref{vpmcf}.

With more precision, throughout this paper we will construct:
\begin{enumerate}
\item Surfaces $\Sigma_i$ ($i = 1,2$) smoothly embedded in $\re^3$ satisfying
\smallskip
\begin{enumerate}
\item[(i)] $\Sigma_1$ is a topological sphere with $H \equiv 0$ on an open set $U \subset \Sigma_1$ (and $H|_{\Sigma_1 \setminus U} > 0)$ (see Theorem \ref{tma-cat})
\smallskip
\item[(ii)] $\Sigma_2$ has genus 1 and its mean curvature  is strictly positive almost everywhere (see Theorem \ref{tma-toro})
\end{enumerate}
which under \eqref{vpmcf} have points with negative mean curvature.
\smallskip
\item A hypersurface $M$ smoothly embedded in $\re^4$ which has non-negative scalar curvature $R \geq 0$, but under \eqref{vpmcf} develops points with $R < 0$ (see Theorem \ref{tma-scal}).
\end{enumerate}

We use a perturbation argument to modify the above examples and obtain a one-parameter family of initial conditions, which are mean convex (have strictly positive mean curvature, respectively) but so that these positivity conditions are lost after evolution under \eqref{vpmcf} with $h(t)$ defined either as in \eqref{oH} or in \eqref{ap-mcf} (see Theorem \ref{tma-pert}).

The proofs to construct all the aforementioned examples include quite easy computations. Hence the relevance of the paper resides in the discovery itself of the counterexamples, and the consequences of their existence in the field: refutation of theorems, widely held {\it beliefs} that are no longer true and the evidence of the need of new (at least different from those in {\sc mcf}) methods to face the singularity analysis of constrained flows.

The paper is organized as follows: Section \ref{pre} contains a miscellanea of preliminary material. In section \ref{bs}, we shall give the examples (1)(i) and (ii). The perturbation argument is carried out in section \ref{smc}. In the last section we prove the analogous result for the scalar curvature.


\section{Preliminaries} \lb{pre}
\subsection{Notation and evolution equations}
Let us first clarify our main notations and sign conventions: for an smooth embedding $F: M \fle \re^{n + 1}$
of an $n$-dimensional manifold as a hypersurface of $\re^{n + 1}$, we consider the principal curvatures $\l_1 \leq \l_2 \leq \cdots \leq \l_n$ which are defined as the eigenvalues of the shape operator $A = (A^i_j)$, where $A^i_j = g^{ik} h_{kj}$, where $g$ is the metric induced on $M$ by the standard Euclidean metric, and $h(X, Y) = \<\nabla_X N,  Y\>$ for $X, Y$ vector fields on $F(M)$, and $\nabla$ is the Levi-Civita connection of $\re^{n + 1}$. Then we define the mean curvature and squared-norm of the second fundamental form as:
$$H = \l_1 + \cdots + \l_n \qquad \text{and} \qquad |A|^2 = \l_1^2 + \cdots + \l_n^2.$$ 
The sign conventions here and in the definition of \eqref{vpmcf} are chosen so that mean convex hypersurfaces (that is, with $H > 0$ everywhere) always move inwards under (unconstrained) {\sc mcf}.

Under \eqref{vpmcf} the above quantities satisfy the following evolution equations (cf.~\cite{Hu87})
\begin{align} \lb{varH}
\parcial{H}{t}  & = \Delta H + |A|^2 (H-h)
  \\ \lb{varA2}
\parcial{|A|^2}{t} &= \Delta|A|^2 - 2 |\nabla A|^2 + 2 |A|^4 - 2 h C,
\end{align}
where $C:= \tr(A^3)$.

 To avoid confusions, we say that a hypersurface is {\it mean convex} if it has positive mean curvature $H > 0$, and we refer to the curvature condition $H \geq 0$ as {\it non-strict mean convexity}.

\subsection{A bit of previous literature}
 The {\sc vp-mcf} was introduced by Gage in \cite{Ga86} (for $n = 1$) and Huisken in \cite{Hu87} for $n\geq 2$. The {\sc ap-mcf} was first studied by Pihan \cite{Pih} (for $n = 1$) and McCoy (\cite{Mc,Mc2} for $n \geq 2$. These papers show that convexity is preserved under the corresponding constrained flow and solutions starting at a convex hypersurface exist for all time and converge to a round sphere. For the evolution of convex hypersurfaces under \eqref{vpmcf} in non-euclidean ambient spaces, see \cite{CaMi1} and for an anisotropic version, see \cite{An}.
 
 A long sought goal is the analysis of singularities for these non-local flows. There are only a few results in this direction \cite{Ath2, Atha-Kan2, CaMi2, CaMi3} and all of them need to work under the extra assumption of rotational symmetry because in this scenario we have some control on the global term.

\subsection{Generalities about hypersurfaces of revolution}
Let $M$ be an embedded  hypersurface of revolution in $\re^{n+1}$ obtained by rotation of a curve $c(u)=(r(u),z(u))$ in the plane $x_1 x_{n+1}$ around the axis $x_{n+1}$. Let  $N = (\dot z(u), - \dot r(u))/|\dot c(u)|$ the unit vector  normal to the curve, and let $\l_1(u)$ be the curvature of $c$ associated to this orientation
\begin{align}\lb{fk1}
\l_1(u) &= \fracc{-1}{|\dot c(u)|^3}\<(\ddot r(u), \ddot z(u)),(\dot z(u),-\dot r(u))\> = \frac{\ddot z \dot r -\ddot r \dot z}{(\dot r^2+ \dot z^2)^{3/2}}.
\end{align} 
The unit normal vector to the hypersurface is obtained by rotation on $N$ and it will be denoted also by $N$. We shall take the orientation of the curve in such a way that $N$ points to the halfspace which do not contains the axis $x_{n+1}$. $\dot c(u)$ is a principal direction of $M$, with principal curvature $\l_1$. All the directions orthogonal to $c'(u)$ (then tangent to the $\mathbb S^n$ spheres which are the parallels of the revolution hypersurface) are also principal, with principal curvature 
\begin{align}\lb{fk2}
\lambda_2(u)=\frac{1}{r(u)} \<N,(1,0)\>= \frac{\dot z(u)}{r(u) \sqrt{\dot r^2+ \dot z^2}} \end{align}
 along $c$ in the plane $x_1x_{n+1}$, and similar formula for the points obtained by rotation of these ones.
 

\section{Loss of non-strict mean convexity} \label{bs}

Let us first highlight that all the examples constructed throughout this paper satisfy $\oH > 0$ and therefore for $h(t)$ defined as in \eqref{ap-mcf}, we have
$$h(t) = \frac{\int_M H^2\, d\mu_t}{\int_M H \, d\mu_t} = \frac{\sfint_M H^2 \, d\mu_t}{\oH},  \qquad \text{where } \quad \fint_M f:= \frac{\int_M f}{|M|}, $$
meaning that $h(t)$ is well defined in all cases, even if our hypersurfaces have some regions where the mean curvature vanishes. Furthermore, we can use  Jensen's inequality to estimate 
\bec \lb{oH_vs_h}
h(t) = \frac{\sfint_M H^2 \, d\mu_t}{\oH} \geq \frac{\(\sfint_M H \, d\mu_t\)^2}{\oH} = \oH.
\eec

\bt \lb{tma-cat}
There exists a surface $\Sigma$ smoothly embedded in $\re^3$ satisfying the following conditions:
\begin{enumerate}
\item $\Sigma$ is a surface of revolution homeomorphic to a sphere.

\item There exists an open subset $U \subset \Sigma$ such that the mean curvature $H$ vanishes on the closure $\overline{U}$ and $H|_{\Sigma \setminus \bar U} > 0$.

\item The solution $\Sigma_t$ of the constrained {\sc mcf} \eqref{vpmcf} (for $h(t)$ defined either as in \eqref{oH} or in \eqref{ap-mcf}) starting at $\Sigma$ develops negative mean curvature at least for a short time, that is, there exists some $\eps > 0$ so that $\inf_{\Sigma} H(\cdot, t) < 0$ for all $t \in (0, \eps)$. 
\end{enumerate}
\et

\bdem
The idea is to construct an embedding $F: \mathbb{S}^2 \hookrightarrow \re^3$ by considering a catenoidal neck to which we paste smoothly two convex caps. We need to perform such a procedure in a way that the resulting sphere has $H\ge 0$; clearly we do not need to worry about the open set $U$ of the catenoidal neck which is not affected by the smooth pasting, since we have $H|_{\overline U} \equiv 0$ 

The catenoid is obtained by rotating a catenary curve $(r(z)= \cosh z,z)$ in $\re^3$. Using the previous formulas \eqref{fk1} and \eqref{fk2} to compute the mean curvature,  the minimality of the catenoid implies, in terms of $r$, that
\bec \lb{min_cat}
\dot r^2 + 1 - \ddot r r = 0,
\eec
where $\dot r$ denotes $\deri{r}{z}$.

Let $a, b \in \re^+$ with $a < b$ and we choose a piece of the catenoid corresponding to $z\in[-b,b]$. To bend the catenoid in the region $z \in [-b, -a) \cup (a, b]$ we introduce a decreasing concave function $\p:[a,b] \flecha \re$ satisfying: 
\bec
\left[\ba{l} \bullet \ \p(a)=1, \p(b)=0, \medskip \\
\bullet \ \p \text{ has all right derivatives vanishing at $z = a$,}\medskip \\ \bullet \ \dot \p, \ddot\p  \text{ are negative on } (a, b), \medskip \\ \bullet \ \p^{-1} \text{ is  concave and }  \frac{d^{s}\p^{-1} }{d r^{s}}\Big|_{r = 0^+} = 0  \text{ with $s = 2k -1$ for all  $k \in \N$.} \ea\right] \lb{phi}
\eec

Now we consider the embedding $F: \mathbb S^2 \hookrightarrow \re^3$ defined by the revolution of the curve $(\rho(z),z)$, where $\rho$ is defined by:
\begin{align}
\rho(z) := \left\{\begin{matrix} r(z) &\text{ for } &-a\le z\le a \smallskip \\ r(z) \p(z) & \text{ for } & a < z \le b \smallskip \\ r(z) \p(-z) & \text{ for } & -b\le z < -a
\end{matrix} \right.. \lb{rho}
\end{align}
It is clear that $\Sigma:=F(\mathbb S^2)$ coincides with the catenoid on the open set $U$ corresponding to $z \in (-a, a)$. We shall obtain $H$ for the other points of $\Sigma$. By symmetry and continuity, it is enough to do the computation for the points where $a< z < b$: 
$$H = \frac{1 + \dot \rho^2 - \ddot \rho \rho}{\rho(1 + \dot \rho^2)^{3/2}}.$$
Therefore to figure out the sign of $H$ in the region corresponding to $z \in (a, b)$ we need to compute
\begin{align*}
1 + \dot \rho^2 - \ddot \rho \rho & = 1 + (\dot r \p + r \dot \p)^2 - (\ddot r \p + 2 \dot r \dot \p + r \ddot \p) r \p
\\ & = 1 + \dot r^2 \p^2 + r^2 \dot \p^2 - r \ddot r \p^2 - r^2 \p \ddot \p
\\ & > \p^2(1 + \dot r^2 - r \ddot r) - r^2 \p \ddot \p >0,
\end{align*}
where for the first inequality we have used that the image of $\p|_{(a, b)}$ lies in $(0, 1)$. The desired positivity  follows by using \eqref{min_cat} and the concavity of $\p$. 

In short, 
$H= 0$ and $|A|^2  > 0$  on $\overline U$  
and $H>0$ outside $\overline U$. This also implies $\Delta H=0$ on $U$ and $h>0$. Substituting all this information in the evolution equation \eqref{varH}, for $t = 0$ we reach
\begin{align*}
\parcial{H}{t}\bigg|_{U}=- |A|^2 h < 0.
\end{align*}
In conclusion, if we flow the embedding $F$ under the {\sc ap-} or {\sc vp-mcf} for each $p\in U$ we can find a small time $\eps(p) >0$ such that $H(p,t)<0$ for all $t\in(0,\eps(p))$. 
\edem

The above theorem gives an example of an embedded sphere with $H\ge 0$ which does not remain with the property $H\ge 0$ after evolving under a constrained {\sc mcf}. But perhaps this phenomenon is not surprising if we take into account that the starting surface had a {\it big chunk} of vanishing mean curvature. For this reason, we also include the following example, where the starting surface has strictly positive mean curvature almost everywhere (more precisely, the initial datum is a torus with elliptic section whose mean curvature vanishes {\it only} along the inner equatorial curve) but still develops negative mean curvature after evolution under \eqref{vpmcf}.

\bt \lb{tma-toro}
We can find a surface $\Sigma$ smoothly embedded in $\re^3$ satisfying the following properties:
\begin{enumerate}
\item $\Sigma$ is a torus of revolution with elliptic section.

\item There exists a curve $\Gamma$ on $\Sigma$ such that $H|_{\Gamma} \equiv 0$ but on $\Sigma \setminus \Gamma$ the mean curvature is strictly positive.

\item Let $\Sigma_t$ be the solution of \eqref{vpmcf} (for $h(t)$ defined either as in \eqref{oH} or in \eqref{ap-mcf}) with $\Sigma_0 = \Sigma$. Then $\Sigma_t$ contains points of negative mean curvature at least for small times: that is, there exists some $\eps > 0$ so that $\inf_{\Sigma} H(\cdot, t) < 0$ for all $t \in (0, \eps)$. 
\end{enumerate}
\et

\bdem
First we consider  the embedding $F: \mathbb S^1\times \mathbb S^1 \hookrightarrow \re^3$ defined by
\be\lb{toro}
F(u,v) = \big((3 + 2 \cos u) \cos v, (3 + 2 \cos u) \sin v, \sqrt{2} \sin u\big).
\ee

 Using the standard formulas from the classical theory of surfaces, we obtain the following expressions for the metric, the outward unit normal vector and the coefficients of the second fundamental form:
$$g_{uu }= 2(1 + \sin^2 u),  \qquad g_{uv}=0, \qquad g_{vv} = ( 3 + 2 \cos u)^2,$$
$$N= \frac{1}{\sqrt{1 + \sin^2 u}} \big(\cos u \cos v, \cos u\sin v, \sqrt{2}\sin u\big),$$
$$h_{uu} = \frac{2}{\sqrt{1 + \sin^2 u}},\qquad h_{uv} =  0, \qquad h_{vv} =  \frac{\cos u(3 + 2 \cos u)}{\sqrt{1 + \sin^2 u}}.$$

We are now in a position to compute the mean curvature:
\begin{align*}
 H & = g^{uu} h_{uu} + g^{vv} h_{vv} = \frac1{(1 + \sin^2 u)^{3/2}} + \frac{\cos u}{\sqrt{1 + \sin^2 u} (3 + 2 \cos u)} 
 \\ & =\frac{\cos^2(u/2) (5+ 2 \cos u - \cos(2 u))}{(3+2 \cos u) (1 + \sin^2 u)^{3/2}}.
\end{align*}
From here we notice that the mean curvature $H$ of $\Sigma:=F(\mathbb S^1 \times \mathbb S^1)$ is a non-negative function of $u$, which only vanishes at $u=\pi$. In other words, $H$ is strictly positive except for the points along the curve $\Gamma$ given by $t \mapsto F(\pi, t)$. 

 Unlike the situation in the proof of Theorem \ref{tma-cat}, here $H$ does not vanish in a neighborhood $U$ of a point, thus we cannot conclude that $\Delta H|_U$ vanish. Rather the contrary, as any point on $\Gamma$ is indeed a minimum for the mean curvature of $\Sigma$, it  holds $\Delta H \ge 0$ in a neighborhood of $\Gamma$ (indeed, we will see that $\Delta H>0$). Looking again at the evolution equation \eqref{varH} and using that $H(\pi) = 0$, we have 
 \bec \lb{varH_pi}
 \parcial{H}{t}\bigg|_{t = 0} (\pi, \cdot) = \Delta H - h \, |A|^2 \Big|_{t=0, u = \pi};
 \eec
hence there is a competition between a positive and a negative term on the right hand side. Here we remark that a careful choice of the coefficients defining $F$ in \eqref{toro} plays a key role to ensure that $\parcial{H}{t}\Big|_{t = 0} (\pi, \cdot) < 0$. The computations to check the latter are carried out in the remaining of the proof (hereafter $t = 0$).

Let us start with the squared norm of the shape operator:
\begin{align} \lb{A2_enpi}
|A|^2\Big|_{u = \pi} = \frac1{(1 + \sin^2 u)^3} + \frac{\cos^2 u}{(1 + \sin^2 u)(3 + 2 \cos u)^2} \bigg|_{u = \pi} = 2. 
\end{align}
On the other hand,
\begin{align*}
\dot H(u) & = \frac{-3 \sin u \cos u}{(1 + \sin^2   u)^{\frac{5}{2}}}
 - \frac{\sin u (1 + \sin^2   u)^{-\frac{1}{2}}}{3 + 2 \cos u} \(1 + \frac{\cos^2 u}{1 + \sin^2   u} - \frac{2 \cos u}{3 + 2 \cos u}\).
\end{align*}
This gives $\dot H(\pi) = 0$. The latter simplifies a lot the computation of the Laplacian, which is basically reduced to the following second derivative (in the computation below we neglect the multiples of $\sin u$, which obviously vanish for $u = \pi$):
$$\ddot H(\pi) = \frac{-3 \cos^2 u}{(1 + \sin^2   u)^{\frac{5}{2}}}  - \frac{\cos u (1 + \sin^2   u)^{-\frac{1}{2}}}{3 + 2 \cos u} \(1 + \frac{\cos^2 u}{1 + \sin^2   u} - \frac{2 \cos u}{3 + 2 \cos u}\)\bigg|_{u = \pi} = 1$$
and from here we get
\begin{align} \lb{lapH_enpi}
\Delta H\Big|_{u = \pi} & = g^{uu} \ddot H(u)\Big|_{u = \pi} = \frac1{2}.
\end{align}

Substituting \eqref{A2_enpi} and \eqref{lapH_enpi} on the right hand side of \eqref{varH_pi}, we reach
\bec \lb{deri_H0}
\parcial{H}{t}\bigg|_{t = 0} (\pi, v) = \frac1{2} - 2 h(0),\eec
whose sign cannot be determined until we compute the specific value of $h(0)$. Here we distinguish the two cases corresponding to the different definitions \eqref{oH} and \eqref{ap-mcf} of $h(t)$:

{\bf Case 1 (VP-MCF):} For $h$ defined as in \eqref{oH} we obtain
$$h(0) =\oH = \frac{\pi}{2 \sqrt{2}} \bigg(\int_0^{\pi/2} \sqrt{1 + \sin^2 \theta} \, d\theta\bigg)^{-1} \geq \frac{\pi}{2 \sqrt{2}} \(\sqrt{2} \,\frac{\pi}{2}\)^{-1} = \frac1{2}.$$
From here, by substitution into \eqref{deri_H0}, we conclude
$$\ds\parcial{H}{t}\bigg|_{t = 0}(\pi,v) \leq -\frac1{2} <0.$$

{\bf Case 2 (AP-MCF):} for $h(t)$ defined as in \eqref{ap-mcf}, thanks to the estimate \eqref{oH_vs_h} we also get 
$$h(0) \geq \oH \geq \frac1{2},$$
and this gives as before $\parcial{H}{t}\big|_{t = 0}(\pi,v) <0$.

In conclusion, both constrained flows satisfy $H\big|_{t = 0}(\pi,v)=0$ and $\parcial{H}{t} \big|_{t = 0}(\pi,v)<0$; therefore, for small $t$, the points on $\Gamma$ evolve to points with $H<0$.
\edem

The following picture illustrates the concrete rate of change of the mean curvature $\parcial{H}{t}\big|_{t=0}$ as a function of $u$ and compares its behaviour under the {\sc vp-mcf} (blue curve) and under the {\sc ap-mcf} (red line):
\begin{center}
\includegraphics[scale=0.6]{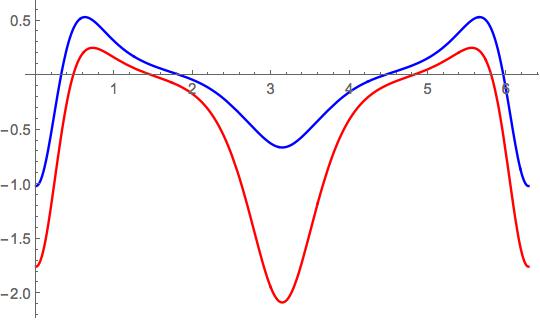} 
\end{center}

 In this section, we have shown explicit examples of surfaces of genus 0 and 1 for which both the {\sc ap-} and {\sc vp-mcf} do not preserve the curvature condition $H \geq 0$. We conjecture that 
 \begin{quote}
{\it there exist similar examples of loss of non-strict mean convexity for surfaces of any genus embedded in $\re^3$.} 
 \end{quote}
The strategy we suggest for the proof follows arguments analogous to those in Theorem \ref{tma-cat}. The starting point would be a minimal surface $\Sigma_g$ of genus $g\ge 0$, with at most a finite number of planar points, and with boundary whose connected components are closed convex curves, near to circles.  The second step would be to paste smoothly a convex cap to each connected component of the boundary of $\Sigma_g$; finally, we would need to check that the resulting surface still has $H \ge 0$.

The following picture illustrates this procedure by taking as $\Sigma_g$ (in this case, $g = 4$) the union of the appropriate number of basic blocks of the Schwarz P minimal surface: 

\includegraphics[width=2.8in]{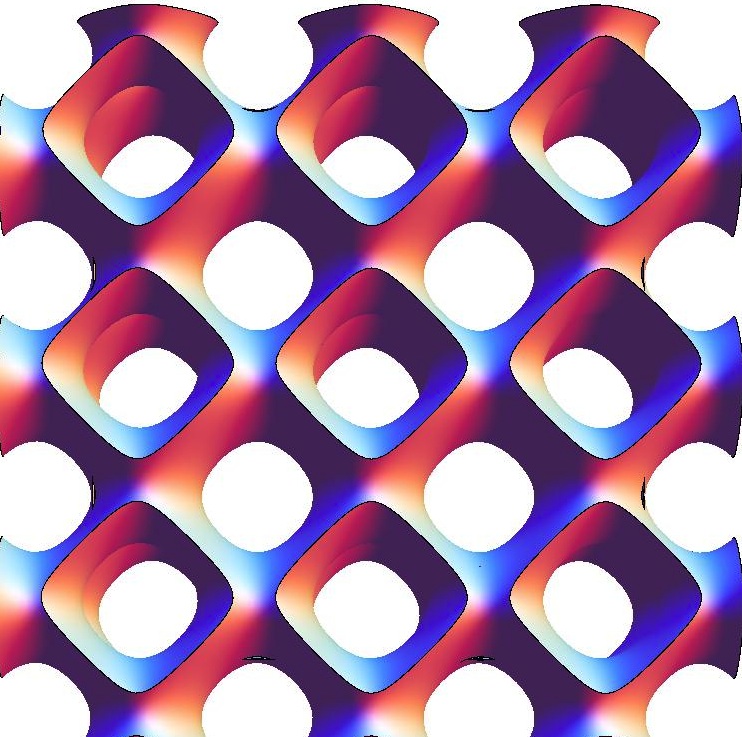}   \includegraphics[width=2.8in]{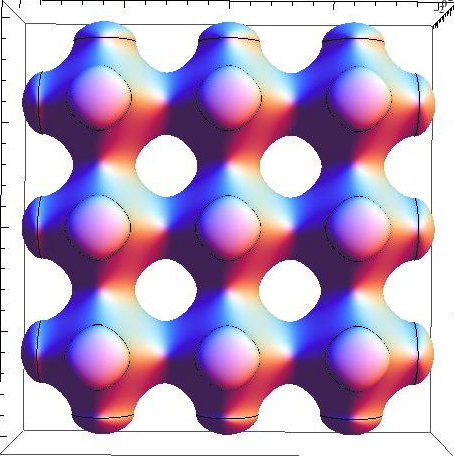}

\noindent The difficulty of our proposal is the following: the lack of an explicit formula for $\Sigma_g$ (in particular, for the curves building its boundary) makes more obscure how to perform the smooth pasting of the convex caps (indeed, even for the picture we used the well-known approximation $\cos x + \cos y + \cos z = 0$ instead of the {\it real} Schwarz P surface). But to overcome such extra complications goes beyond the scope of the present paper, since we have enough with either Theorem \ref{tma-cat} or Theorem \ref{tma-toro} to construct our desired examples of loss of mean convexity under \eqref{vpmcf}. This is actually the content of the next section.

\section{A perturbation argument giving loss of mean convexity}\label{smc}

The goal is to produce examples of embedded spheres with $H>0$ initially that also evolve by {\sc ap-} or {\sc vp-mcf} giving points with $H<0$. The idea is to perturb the examples $\Sigma$ constructed in the previous section by flowing them under the (unconstrained) {\sc mcf} for very small time $s \in (0, \eps)$. This will give a family of initial surfaces $\Sigma^s$ with $\Sigma^0 = \Sigma$ and  $H_{\Sigma^s} > 0$ for $s > 0$. These new starting data will still develop points of negative mean curvature under \eqref{vpmcf} by using the smooth dependence of the initial condition for quasilinear parabolic equations proved by C.~Mantegazza and L.~Martinazzi in \cite{MaMa}.

\bt \lb{tma-pert}
 There exists a positive time $t>0$ such that for some $\sigma > 0$ we can find a one-parameter family of mean convex surfaces $\{F^s: \mc S \hookrightarrow \re^3\}_{s \in (0, \sigma)}$  smoothly embedded in $\re^3$ and satisfying that, under the constrained {\sc mcf} \eqref{vpmcf} (with $h(t)$ defined either as in \eqref{oH} or in \eqref{ap-mcf}), it holds:
 $$\inf_{\mc S} H(F^s(\cdot, t)) < 0.$$
\et

\bdem
Let $F:\mathcal S \hookrightarrow \re^3$ be either the embedding  defined in the proof of Theorem \ref{tma-cat} or the one from Theorem \ref{tma-toro}. Recapitulating, if we take $\{F_t:=F(\cdot,t)\}_{t\in [0,T)}$ the solution of \eqref{vpmcf} with initial condition $F_0 = F$, we already know that
\begin{enumerate}
\item[(i)] $H(\cdot,0) \geq 0$, and

\item[(ii)] there is some time $\tau > 0$ such that $\ds\inf_{\mc S}H(\cdot , t) < 0$  for every  $t\in (0,\tau).$
\end{enumerate}
The above property (ii) implies that for any  interval $[\tau_1,\tau_2]\subset (0,\tau)$, there is a constant $c >0$ such that 
\bec \lb{last_S}
\inf_{\mc S}H(\cdot, t) < -c  \qquad \text{ for all} \quad t\in [\tau_1,\tau_2].
\eec

On the other hand, let us consider  $X: \mc S \times [0, S) \fle \re^3$ a solution of the (uncostrained) {\sc mcf}:
\bec \lb{mcf}
\left\{\ba{ll} \ds\parcial{X}{s} (p, s) = - H N(p, s), &\ p\in \mc S, \ s \in (0, S) \medskip \\ X(\cdot, 0) = F_0\ea \right. .
\eec
 Taking into account the property (i) above, by application of the strong maximum principle to the evolution equation \eqref{varH} with $h =0$, we get that the mean curvature $H_X$ of $X(\cdot,s)$ satisfies
\bec \lb{mean_co}
H_{X}(\cdot,s) >0  \qquad \text{ for any } \qquad s \in (0, S).
\eec

Next we look at  $F^s_t:= F^s(\cdot,t)$, where $F^s: \mc S \times [0, T_\ast) \fle \re^3$ is the solution of \eqref{vpmcf} with $F^s_0 = X(\cdot,s)$. For some $t \in [\tau_1, \tau_2]$ fixed,  the continuity of  $\inf_{\mc S} H$ (regarded as a map on the space of all the immersions of $\mc S$ into $\re^3$) implies that, given $\tilde \eps = c/2 > 0$ (where $c$ is the constant coming from \eqref{last_S}), we can find  $\tilde \delta = \tilde \delta(c, t)> 0$ such that the following happens: If 
\bec \lb{todo}
d_{C^2}(F_t, F_t^s) < \tilde{\delta} \qquad \text{ for all } \quad s \in (0, \sigma),
\eec
then we conclude
\bec \lb{final_concl}
\inf_{\mc S}H(F^s(\cdot,t) )< -c/2 .
\eec

It remains to check \eqref{todo}. With this goal, we recall that by \cite[Theorem 1.1]{MaMa}, the solution of \eqref{vpmcf} depends continuously (in the $C^\infty$-topology) on the smooth initial condition. Accordingly, given $k\in \N$ and $\eps>0$ arbitrary, there exists a $\delta = \delta(k,\eps)>0$ satisfying 
\bec \lb{MaMa}
d_{C^k}(F_t,F^s_t)<\eps \quad \text{for all } \quad t\in[0, \min\{T, T_\ast\}), \quad \text{whenever} \quad d_{C^k}(F_0,F_0^s)<\delta.
\eec

Hereafter we fix $k = 2$ and choose $\eps = \tilde \delta > 0$. Then we can find $\delta = \delta(\eps) > 0$ satisfying the property in \eqref{MaMa}. By the smooth dependence with respect of the time parameter $s$ of the map $X: \mc S \times [0, S) \fle \re^3$ solving \eqref{mcf}, we can find a small time $\sigma = \sigma(\delta) > 0$ such that
$$d_{C^2}(F_0, F_0^s) < \delta \qquad \text{for all } \quad s \in [0, \sigma).$$ 
By application of \eqref{MaMa}, we obtain the inequality in \eqref{todo}, which yields the conclusion \eqref{final_concl}. From here the statement follows, since the embedded surfaces $\{F^s_0(\mc S)\}_{s \in (0, \sigma)}$ are mean convex (recall \eqref{mean_co}) but such a property is lost under \eqref{vpmcf} for some time $t >0$.
\edem




\section{Non-preservation of positivity of scalar curvature for hypersurfaces in $\re^4$}\lb{Smo}

After the failure of the preservation of the mean convexity under {\sc ap-} and {\sc vp-mcf}, we wonder if the positivity of the scalar curvature ($R > 0$) is preserved or not under these constrained flows. This is a natural question by e.g.~the following two reasons:

(1) Unlike mean convexity, $R > 0$ is preserved for $n = 2$. In fact, for surfaces in $\re^3$ the preservation of $R = \l_1 \l_2 > 0$ is equivalent to the preservation of  convexity, and the latter was proved by G.~Huisken in \cite{Hu87}. 

(2) $R > 0$ corresponds to the positivity of the second symmetric polynomial $\sigma_2 :=\sum_{i < j} \l_i \l_j = \frac{R}{2}$ (notice that also $H = \sigma_1$). Therefore, as we should plan to study a curvature condition which is stronger than $H > 0$ (since we already know that this is a non-preserved condition under \eqref{vpmcf}), the next natural candidate is $\sigma_2 > 0$.

As suggested by the title of the section, we provide a negative answer: we show an example of a revolution hypersurface in $\re^4$ for which the preservation of the positivity of the scalar curvature fails again. The ideas for the construction of this example are similar to those in the proofs of Theorem \ref{tma-cat} and Theorem \ref{tma-pert}. 

Roughly speaking, here we consider a paraboloidal (instead of catenoidal) neck defined by rotation of the curve
\bec \lb{parab}
(r(z), z) \qquad \text{with} \qquad r(z) = 2 + z^2/8
\eec
around the $x_4$ axis in $\re^4$ and the central cross section is in this case a 2-sphere of radius 2. In other words, we depart from a revolution hypersurface $M$ in $\re^4$ generated by the curve defined in \eqref{parab} or, in parametric coordinates, by $c(u) = (r(u), z(u))$ with
\bec \lb{def_rz}
r(u) = 2 \cosh^2 u \qquad \text{ and } \qquad z(u)=4 \sinh u.
\eec 
Applying \eqref{fk1} and \eqref{fk2}, the principal curvatures of $M$ are given by
\begin{align} \lb{pcurv}
\l_1= -\l(u)  \quad \text{ and } \quad  \l_2 = \l_3 = 2 \l(u), \qquad \text{with} \quad \l(u) = \frac1{4 \cosh^3u}.
\end{align}
From here, the scalar curvature of  $M$ is 
\bec \lb{Rnull}
R = 2(2 \l_1 \l_2 + \l_2^2) = 2 \l_2(2 \l_1 + \l_2) =0. 
\eec
Using again \eqref{fk1}, \eqref{fk2} and that $\l_2 >0$, we can rephrase \eqref{Rnull} by saying that 
\bec \lb{Sr}
S(r):= 2 r(\dot r \ddot z - \ddot r \dot z) + \dot z(\dot r^2 + \dot z^2) =0.
\eec

By adding two caps to $M$, we will construct in the next theorem  a closed 3-manifold $M^\ast$ with $R \geq 0$ but that under \eqref{vpmcf} develops points of negative scalar curvature.
\bt \lb{tma-scal}
There exist $t'>0$,  $\sigma >0$ and  a 1-parameter family of smooth embeddings $F^s: \mathbb S^3 \hookrightarrow \re^4$, with $s \in [0, \sigma)$ such that:

$\circ$ For all $s>0$ the hypersurfaces $M^s:= F^s(\mathbb S^3)$ have strictly positive scalar curvature.

$\circ$ The solution $M^s_t$ of \eqref{vpmcf} (with $h(t)$ defined either as in \eqref{oH} or in \eqref{ap-mcf}) starting at $M^s$ satisfies $\inf_{M^s} R(\cdot, t') < 0$.
\et

\bdem
Our starting point is a revolution hypersurface $M$ in $\re^4$ whose generating curve $c(u) = (r(u), z(u))$ is defined as in \eqref{def_rz}. We  choose a piece of $M$ corresponding to $z(u)\in[-b,b]$, for some $b \in \re^+$.

Now we modify $M$ in the following way: we take $a \in \re^+$ with $a < b$ and consider the embedding $F: \mathbb S^3 \hookrightarrow \re^4$ defined by the revolution of a curve $(\rho(u),z(u))$, with
$$\rho(u):= \left\{\ba{lcl} r(u)\, \p(-z(u)), & \text{if} & z(u) \in [-b, -a) \smallskip \\
r(u) & \text{if} & z(u) \in [-a, a] \smallskip \\
r(u)\,\p(z(u)) & \text{if} & z(u) \in (a, b]\ea\right.$$
and $\p:[a,b] \fle \re$ is the bending function satisfying the properties in \eqref{phi}. Clearly, $M^\ast:= F(\mathbb S^3)$ coincides with $M$ on the open set $U$ where $-a<z(u)<a$. 

The scalar curvature $R$ for the  points of $M^\ast$ corresponding to $z(u) \in (a, b)$ is then 
$$R = 2 \l_2(2 \l_1 + \l_2) = 2 \l_2 \frac{S(\rho)}{\rho (\dot \rho^2 + \dot z^2)^{3/2}},$$
where $S(\rho)$ is given by formula \eqref{Sr} (writing $\rho, \dot \rho, \ddot \rho$ instead of $r, \dot r, \ddot r$). As $\l_2 > 0$, we need to study the sign of 
\begin{align*}
S(\rho) & = 2 \ddot z (\dot\p \,\dot z \, r + \p \,\dot r) \, r\, \p - 2 \dot z\, r \, \p (\ddot \p \, \dot z^2 r + \dot \p \ddot z \, r + 2 \dot \p\, \dot r \, \dot z + \p \, \ddot r) + \dot z(\dot \p\, \dot z \, r + \p \, \dot r)^2 + \dot z^3 
\\ & > \p^2 S(r) - 2 \p \, \ddot \p \, r^2 \dot z^3 + r^2 \dot \p^2 \dot z^3 - 2 \p \, \dot \p\, \dot z^2 \, r \, \dot r \geq 0,
\end{align*}
where the first inequality follows by using $0 < \p|_{(a,b)} < 1$ and for the non-negativity we have applied  \eqref{Sr}, the concavity of $\p$ and that $\dot \p < 0$.

In summary, we have constructed $M^\ast$ a closed hypersurface of $\re^4$ with non-negative scalar curvature and such that it has an open subset $U$ where $R$ vanishes identically. 

Next, we are interested in the rate of change of $R$ under \eqref{vpmcf}. With the goal of computing the corresponding evolution formula, notice that, by the Gauss equation, the scalar curvature of a hypersurface in $\re^{n + 1}$ is given by 
$$R =  2 \sum_{i< j} \l_i\, \l_j = H^2-|A|^2.$$
Taking into account the evolution equations \eqref{varH} and \eqref{varA2}, we obtain
\begin{align}\lb{pstV}
\parcial{R}{t} &=  \Delta R + 2 (|\nabla A|^2 -|\nabla H|^2)+  2\ |A|^2 R - 2 h (H |A|^2 - C).
\end{align}

Hereafter we compute the terms of the right hand side of the above formula for $z(u) \in (-a, a)$. Using \eqref{pcurv}, we can easily compute the last term:
\begin{align}
H |A|^2  - C &= (-\l+4 \l) (\l^2+2 (2\l)^2) - (2 (2 \l)^3 - \l^3) =  12 \l^3.  \lb{exV}
\end{align} 
In order to compute the gradient terms of \eqref{pstV}, we introduce the coordinates $(u,\theta)\in \re\times \mathbb S^2$ to parametrize $M$ so that the induced metric of $M$ is given by
$$g_M = 4 \cosh^4 u (4\ du^2 + g_{_{\mathbb S^2}}) = r^2(u) (4\ du^2 + g_{_{\mathbb S^2}}),$$ 
where $g_{\mathbb S^2}$ is the standard metric on the unit sphere $\mathbb S^2$. This allows us to use the well-known formula (see e.g.~\cite[page 90]{Pet}) for the covariant derivatives of a metric conformal to the product metric $4\ du^2 + g_{\mathbb S^2}$ of $\re\times \mathbb S^2$. Indeed,
using that $\nabla_{\partial_u} \partial_u$ is proportional to $\partial_u$ and $A \partial_u = - \lambda \partial_u$ we get
\begin{align} \lb{Auu}
(\nabla_{\partial_u} A)\partial_u &= \nabla_{\partial_u}(A\partial_u)- A(\nabla_{\partial_u}\partial_u) = -\dot\l \,\partial_u -\l \nabla_{\partial_u} \partial_u + \l \nabla_{\partial_u} \partial_u = - \dot\l \, \partial_u. 
\end{align}
For any vector field $w$ tangent to $\mathbb S^2$,  after realizing that $\nabla_{\partial_{u}} w = \nabla_w \partial_u \frac{\dot r}{r} w$ and $Aw = 2 \l w$, we deduce
\begin{align} \lb{Aui}
(\nabla_{\partial_u} A)w &= \partial_u(2 \l) w + 2 \l \nabla_{\partial_u} w - 2 \l \nabla_{\partial_u} w =  2 \dot \l \, w ;
\end{align}
\begin{align} \lb{Aiu}
(\nabla_w A) {\partial_u}& =  - \l \nabla_w{\partial_u} - 2 \l \nabla_w{\partial_u} = -3 \l \nabla_w{\partial_u} = - 3 \l \, \frac{\dot r}{r}  \, w.
\end{align}
If we notice that for any local frame $\{v, w\}$ of the unit sphere $\mathbb{S}^2$ we get
$$\nabla_v w = \nabla^{\mathbb S}_v w - \frac{\dot r}{4 r} g_{_{\mathbb S^2}}(v, w) \, \partial_u,$$
where $\nabla^{\mathbb S}$ denotes the Levi-Civita connection of $\mathbb S^2$ with its standard metric, then we have
\begin{align} \lb{Aii}
(\nabla_v A) w&= 2 \l \nabla_v w - A(\nabla^{\mathbb S}_v w) + \frac{\dot r}{4 r} g_{_{\mathbb S^2}}(v, w) A(\partial_u) \nn
\\ & = 2 \l \nabla^{\mathbb S}_v w  - \frac{\l\, \dot r}{2 r}  g_{_{\mathbb S^2}}(v, w) \, \partial_u - 2 \l \nabla^{\mathbb S}_v w - \frac{\l \, \dot r}{4 r} g_{_{\mathbb S^2}}(v, w) \partial_u \nn
\\ & = - \frac{3 \l \,\dot r}{4 r}  g_{_{\mathbb S^2}}(v, w) \, \partial_u.
\end{align}

Next we choose a local frame  $\{e_1=\partial_u,  e_2, e_3\}$, where $\{e_2,e_3\}$ is a local orthonormal frame of the unit sphere $\mathbb S^2$. Using the computations \eqref{Auu}, \eqref{Aui}, \eqref{Aiu}, \eqref{Aii}, the symmetries of $\nabla A$ and that $g^{22} = g^{33} = 4 g^{uu}$, we can write
\begin{align}
|\nabla A|^2 &= g^{ij} g^{k\ell}  \<(\nabla_{e_i} A)e_k, (\nabla_{e_j} A)e_\ell\>   \nn \\
&= (g^{uu})^2 |(\nabla_{\partial_u} A)\partial_u|^2 + 2 g^{uu} g^{22} \(|(\nabla_{\partial_u} A)e_2|^2 + |(\nabla_{e_2} A)\partial_u|^2\)  
  + 2 (g^{22})^2 |(\nabla_{e_2} A)e_2|^2\nn \\
&= g^{uu} \dot \l^2 + 2 g^{uu}\Big[(2 \dot\l)^2 +  (3 \l \dot r/r)^2\Big]  + 2 (g^{22})^2 \Big(\frac{3 \l\, \dot r}{4 r}\Big)^2 g_{uu} \nn \\
& = 18 g^{uu} \Big(2 \l^2 \frac{\dot r^2}{r^2} + \frac{\dot \l^2}{2}\Big).
\end{align}
After checking that $|\nabla H|^2 =  g^{uu} (3\dot \l)^2$, we conclude
\bec \lb{gradT}
|\nabla A|^2 -|\nabla H|^2  = 9 g^{uu} \Big(4 \l^2 \frac{\dot r^2}{r^2} + \dot \l^2 - \dot \l^2\Big) = 36 g^{uu} \frac{\l^2 \dot r^2}{r^2} =  \Big(3 \frac{\l \, \dot r}{r^2}\Big)^2.
\eec
Now we substitute \eqref{Rnull}, \eqref{exV} and \eqref{gradT} on the right hand side of \eqref{pstV}; this yields
\begin{align}\lb{menoro}
\parcial{R}{t}=  2 \Big(3 \frac{\l \, \dot r}{r^2}\Big)^2 - 24\, h(t) \l^3.
\end{align}
Here we have used that $\Delta R$ vanishes on the open set $U$ where $z(u) \in (-a, a)$  (since $R|_{\overline U} \equiv 0$). Notice that the first term  also vanishes for $u = 0$ because it is a multiple of $\dot r(u)$ and, therefore,   includes powers of $\sinh u$.  On the other hand, by \eqref{pcurv} we also have $H|_{\overline U} > 0$ and (as one can easily check that $H \geq 0$ on $M^\ast$) hence $\oH > 0$. The latter implies $h(t) > 0$ (for $h$ as in \eqref{oH} this follows by definition and, if we take $h$ as in \eqref{ap-mcf}, because of the estimate \eqref{oH_vs_h}).

Taking into account the above considerations, we substitute \eqref{pcurv} into \eqref{menoro} to deduce
\bec \lb{dR_u0}
\parcial{R}{t}\bigg|_{u = 0} = - \frac3{8} h(t) \frac1{\cosh^9 u} \bigg|_{u = 0} = - \frac3{8} h(t) < 0.
\eec
This indicates that under \eqref{vpmcf} we can find points of $M^\ast$ in a neighborhood of $u = 0$ which develop negative scalar curvature for small positive times.

Finally, the statement follows by using the same perturbation argument as in the proof of Theorem \ref{tma-pert}. The only modification is that here we need to use that the solution of (unconstrained) {\sc mcf} starting at a hypersurface with $R \geq 0$ has positive scalar curvature for any $t > 0$ (cf.~\cite[Theorem 3]{HuWu}).

\edem

{\bf \small Acknowledgments:} Work partially supported  by DGI (Spain) and FEDER Project   MTM2010-15444, MTM2013-46961-P
  and the Generalitat Valenciana Project  PROMETEOII/2014/064.  
	
	We wish to thank Burkhard Wilking for helpful discussions about this topic.

\bibliographystyle{alpha}

\end{document}